\documentclass{article}%
\usepackage{amsmath}
\usepackage{amsfonts}
\usepackage{amssymb}
\usepackage{graphicx}%
\providecommand{\U}[1]{\protect\rule{.1in}{.1in}}
\newtheorem{theorem}{Theorem}

\newtheorem{definition}[theorem]{Definition}
\newtheorem{example}[theorem]{Example}

\begin{document}

\title{Notes on Degenerate Curves in Pseudo-Euclidean Spaces of Index Two}
\author{Mehmet G\"{o}\c{c}men, Sad\i k Kele\c{s}} \maketitle

\begin{abstract}
In this paper we deal with curves with degeneration degree two in
pseudo-Euclidean spaces of index two. We characterize Bertrand curves. We show
a correspondence between the evolute of a null curve and the involute of a
certain spacelike curve in the $6-$dimensional pseudo-Euclidean space of index
two. Also we characterize pseudo-spherical null curves in the $n-$dimensional
pseudo-Euclidean space of index two in terms of the curvature functions.

\end{abstract}

\noindent\textbf{Mathematics Subject Classification:} 53A04, 53B30.

\noindent\textbf{Keywords and phrases:} Degeneration degree, Radical
of a Curve, Index of a Curve, Index Sequence, Nullity Degree
Sequence, Bertrand curve.

\section{Introduction}

Due to the growing importance of degenerate geometry in mathematical
physics ( null curves and null hypersurfaces, etc ), mathematicians
has been trying to get good Frenet frames for degenerate curves in
pseudo-Euclidean spaces. The geometry of null hypersurfaces in
space-times has played an important role in the development of
general relativity, as well as in mathematics and physics of
gravitation. It is necessary, e.g. to understand the causal
structure of space-times, black holes, assymptotically flat systems
and gravitational waves. An initial point to study null surfaces, or
in general null hypersurfaces, consists of investigating the curves
that lie in those hypersurfaces. In this sense null curves in
Lorentzian space forms have been studied by several authors
(\cite{bejancu}, \cite{bonnor}). The pseudo-Euclidean spaces of
index two were studied by Duggal and Jin \cite{jin}. Ferrandez,
Gimenez and Lucas introduced Cartan reference along a degenerate
curve. They obtained several different types of degenerate curves
and present existence, uniqueness and congruence theorems
\cite{intex-2}. Sakaki \cite{sakaki} showed a correspondence between
the evolute of a null curve and the involute of a certain spacelike
curve in the $4-$dimensional Minkowski space. He also characterized
pseudo-spherical null curves in the $n-$dimensional Minkowski space
in terms of the curvature functions.

In this paper we discuss null curves in $R_{2}^{n}.$ We characterize Bertrand
curves in $R_{2}^{5}$ and pseudo-spherical null curves in $R_{2}^{n}$ in terms
of the curvature functions (see \cite{ceylan}, \cite{duggal}, \cite{matsuda}).
Then we define the evolute of a null curve in $R_{2}^{6}$ and the involute of
a spacelike curve in $R_{2}^{6},$ and show a correspondence between them which
is similar to that between the plane evolute and involute.

\section{Preliminaries}

Let $V$ be an $n-$ dimensional real vector space endowed with a symmetric
bilinear mapping $g:V\times V\rightarrow R.$ We will say that $g$ is
degenerate on $V$ if there exists a vector $\xi\neq0$ of $V$ such that
\[
g\left(  \xi,v\right)  =0,\text{ \ \ for all }v\in V.
\]
Otherwise, $g$ is said to be non-degenerate. The radical of $\left(
V,g\right)  $ is the subspace of $V$ defined by
\[
Rad\text{ }V=\left\{  \xi\in V;g\left(  \xi,v\right)  =0\text{ for all }v\in
V\right\}  .
\]
It is clear that $V$ is non-degenerate if and only if $Rad$
$V=\left\{ 0\right\}  .$ A Pseudo-Euclidean space $\left( V,g\right)
$ will be an $n-$dimensional real vector space $V$ equipped with a
symmetric non-degenerate bilinear mapping $g$. The dimension $q$ of
the largest subspace $W\subset V$ on which $g\left\vert _{W}\right.
$ is definite negative is called the index of $g$ on $V$. $\left(
V,g\right)  $ will be denoted by $R_{q}^{n}$.

Let $B=\left\{  V_{1},...,V_{n}\right\}  $ be an ordered basis of a
pseudo-Euclidean space and let $r_{i}$ and $q_{i}$ be the dimension of the
radical and the index of \emph{span} $\left\{  V_{1},...,V_{i}\right\}  $ for
all $i,$ respectively. The sequence $\left\{  r_{i};0\leqslant i\leqslant
n\right\}  $ and $\left\{  q_{i};0\leqslant i\leqslant n\right\}  ,$ where
$r_{0}=q_{0}=0,$ will be called the \emph{nullity degree sequence} and the
\emph{index sequence} of the basis $B.$ It is easy to see that $\left\vert
r_{i}-r_{i-1}\right\vert $ and $q_{i}-q_{i-1}$ are either $0$ or $1$, for all
$i=1,2,...,n,$ as well as $r_{n}=0$ and $q_{n}=q.$

Let $B=\left\{  V_{1},...,V_{n}\right\}  $ be an ordered basis of
pseudo-Euclidean space and let $\left\{  r_{i};1\leqslant i\leqslant
n\right\}  $ be the nullity degree sequence. The positive number
\[
r=\frac{1}{2}%
{\displaystyle\sum\limits_{i=1}^{n}}
\left\vert r_{i}-r_{i-1}\right\vert
\]
is said to be the \emph{degeneration degree}\textbf{ }of the basis $B$. Let
$R_{2}^{n}$ be a pseudo- Euclidean space of index two and let $\alpha
:I\rightarrow R_{2}^{n}$\ be a differentiable curve in $R_{2}^{n}$. Assume
that $A=\left\{  \alpha^{\prime}\left(  t\right)  ,...,\alpha^{\left(
n\right)  }\left(  t\right)  \right\}  $ is a linearly independent system for
all $t\in I$, and, for all $i,$ $r_{i}\left(  t\right)  $ and $q_{i}\left(
t\right)  $ are constant for all $t\in I,$ where $\left\{  r_{i}\left(
t\right)  ;0\leqslant i\leqslant n\right\}  $ and $\left\{  q_{i}\left(
t\right)  ;0\leqslant i\leqslant n\right\}  $ stand for the nullity degree and
the index sequences of the basis $A$. In this case these sequences will be
called nullity degree and the index sequences of the curve $\alpha,$
respectively, and the degeneration degree $r$ ($=$ constant) of $A$ will be
called the degeneration degree of the curve $\alpha.$

With the above notation, a curve $\alpha:I\rightarrow R_{2}^{n}$ is said to be
a \emph{degenerate curve} if $r>0$ \cite{intex-2}.

To study curves of $R_{2}^{n}$, we use the metric $\left\langle ,\right\rangle
$ defined as follows%
\[
\left\langle x,y\right\rangle =-x^{1}y^{1}-x^{2}x^{2}+x^{3}x^{3}+x^{4}%
x^{4}+...+x^{n}x^{n}%
\]

\noindent for all vectors $x,y\in R_{2}^{n};$ $x=\left(  x^{1},x^{2}%
,x^{3},...,x^{n}\right)  ,$ $y=\left(  y^{1},y^{2},y^{3},...,y^{n}\right)  ,$
$x^{i},y^{i}\in R,$ $1\leqslant i\leqslant n.$

\section{Bertrand Curves in $R_{2}^{5}$}

In the whole section we restrict our search to the family-type of
curves whose nullity degree sequence is $\left\{
0,1,2,2,1,0\right\}  $ to study Bertrand curves in R$_{2}^{5}$,$.$
Now $\alpha$ be a null Cartan curve with the degeneration degree two
in $R_{2}^{5}.$ We shall determine under which conditions this curve
is a Bertrand curve. For $\alpha$ assume that $\left\{
\alpha^{\prime},\alpha^{\prime\prime},\alpha^{\left(  3\right)  }%
,\alpha^{\left(  4\right)  },\alpha^{\left(  5\right)  }\right\}  $ is
positively oriented, the sets $\left\{  \alpha^{\prime},\alpha^{\prime\prime
},\alpha^{\left(  3\right)  },\alpha^{\left(  4\right)  },\alpha^{\left(
5\right)  }\right\}  $ and $\left\{  L_{1},L_{2},W_{3},N_{2},N_{1}\right\}  $
have the same orientation. We choose $\left\langle L_{1},N_{1}\right\rangle
=1,$ $\left\langle L_{2},N_{2}\right\rangle =-1.$ Then the Cartan equations
are as follows:%
\begin{align*}
\alpha^{\prime}  &  =L_{1},\text{ \ }L_{1}^{\prime}=L_{2},\text{ \ }%
L_{2}^{\prime}=W_{3},\text{ \ }W_{3}^{\prime}=-k_{1}L_{2}+N_{2},\\
N_{2}^{\prime}  &  =k_{2}L_{1}+N_{1}-k_{1}W_{3},\text{ \ }N_{1}^{\prime}%
=k_{2}L_{2}.
\end{align*}

\begin{definition}
Let $\left(  \alpha,\overline{\alpha}\right)  $ be a pair of \ framed null
Cartan curves in $R_{2}^{5}$, with pseudo-arc parameters $s$ and $\overline
{s},$ respectively. This pair is said to be a null Bertrand pair if their
spacelike vectors $W_{3}$ and $\overline{W}_{3}$ are linearly dependent. The
curve $\overline{\alpha}$ is called a Bertrand mate of $\alpha$ and vice
versa. A framed null curve is said to be a null Bertrand curve if it admits a
Bertrand mate. To be precise, a null Cartan curve $\alpha$ in $R_{2}^{5}$
$\left(  \alpha:I\rightarrow R_{2}^{5}\right)  $ is called a Bertrand curve if
there exist a null Cartan curve $\overline{\alpha}$ $\left(  \overline{\alpha
}:\overline{I}\rightarrow R_{2}^{5}\right)  ,$ distinct from $\alpha,$ and a
regular map $\varphi:I\rightarrow\overline{I}\left(  \overline{s}%
=\varphi\left(  s\right)  ,\frac{d\varphi\left(  s\right)  }{ds}\neq0\text{
for all }s\in I\right)  $ such that the spacelike vectors $W_{3}$ of $\alpha$
and $\overline{W}_{3}$ of $\overline{\alpha}$ are linearly dependent at each
pair of corresponding points $\alpha\left(  s\right)  $ and $\overline{\alpha
}\left(  \overline{s}\right)  =\overline{\alpha}\left(  \varphi\left(
s\right)  \right)  $ under $\varphi$ (\cite{duggal}, \cite{matsuda}).
\end{definition}

\begin{theorem}
Let $\alpha$ be a null Cartan curve in $R_{2}^{5}$ having the
degeneration degree two and nullity degree sequence $\left\{
0,1,2,2,1,0\right\}  $. Then
$\alpha$ is a Bertrand curve iff the curvature functions of $\alpha$ satisfy%
\[
k_{1}=k_{2}=0.
\]

\end{theorem}

\noindent\textbf{Proof}. Assume that $\alpha$ is a Bertrand curve. Then there
exist a Bertrand mate $\overline{\alpha}$ of $\alpha$ and the normal lines of
vector fields $W_{3}$ and $\overline{W_{3}}$ of $\alpha$ and $\overline
{\alpha}$, respectively coincide. So we set%
\begin{equation}
\overline{\alpha}\left(  \overline{s}\right)  =\alpha\left(  s\right)
+\mu\left(  s\right)  W_{3}\left(  s\right)  \label{a1}%
\end{equation}

\noindent If we differentiate (\ref{a1}) with respect to $s$, we get
\begin{equation}
\overline{L_{1}}\left(  \overline{s}\right)  \frac{d\overline{s}}{ds}%
=L_{1}\left(  s\right)  +\mu^{\prime}\left(  s\right)  W_{3}\left(  s\right)
-\mu\left(  s\right)  k_{1}\left(  s\right)  L_{2}\left(  s\right)
+\mu\left(  s\right)  N_{2}\left(  s\right)  \label{a2}%
\end{equation}

\noindent Since $\left\langle W_{3}\left(  s\right)  ,\overline{L_{1}}\left(
\overline{s}\right)  \right\rangle =0,$ we have $\mu^{\prime}=0.$ So $\mu$ is
a nonzero constant number. In this case, (\ref{a2}) becomes%
\begin{equation}
\overline{L_{1}}\left(  \overline{s}\right)  \frac{d\overline{s}}{ds}%
=L_{1}\left(  s\right)  -\mu\left(  s\right)  k_{1}\left(  s\right)
L_{2}\left(  s\right)  +\mu\left(  s\right)  N_{2}\left(  s\right)  \label{a3}%
\end{equation}

\noindent If we apply the metric $\left\langle ,\right\rangle $ on each side
of the equation (\ref{a3}), we get
\[
2\mu^{2}k_{1}=0,
\]

\noindent so we obtain%
\[
k_{1}=0.
\]

\noindent If we differentiate (\ref{a3}) once more, we have
\begin{equation}
\overline{L_{2}}\left(  \overline{s}\right)  \left(  \frac{d\overline{s}}%
{ds}\right)  ^{2}+\overline{L_{1}}\left(  \overline{s}\right)  \frac
{d^{2}\overline{s}}{ds^{2}}=L_{2}\left(  s\right)  +\mu\left[  k_{2}\left(
s\right)  L_{1}\left(  s\right)  +N_{1}\left(  s\right)  \right]  . \label{a4}%
\end{equation}

\noindent Applying the metric $\left\langle ,\right\rangle $ on each side of
the equation (\ref{a4}), we get
\[
2\mu^{2}k_{2}=0.
\]
So we obtain%
\[
k_{2}=0.
\]

Conversely, let $\alpha$ be a null Cartan curve in $R_{2}^{5}$ whose cutvature
functions satisfy%
\[
k_{2}=k_{1}=0.
\]

\noindent Then we will show that such a curve is Bertrand curve. Assume that a
curve $\overline{\alpha}$ in $R_{2}^{5}$ defined as follows%
\begin{equation}
\overline{\alpha}\left(  s\right)  =\alpha\left(  s\right)  +\mu W_{3}\left(
s\right)  \label{a5}%
\end{equation}

\noindent$\mu$ is a nonzero constant. If we differentiate (\ref{a5}) with
respect to $s$ and use the Frenet equations, we obtain
\[
\frac{d\overline{\alpha}\left(  s\right)  }{ds}=L_{1}\left(  s\right)  +\mu
N_{2}\left(  s\right)
\]
It can be seen that $\frac{d\overline{\alpha}\left(  s\right)  }{ds}\neq0.$
Because if $\frac{d\overline{\alpha}\left(  s\right)  }{ds}=L_{1}\left(
s\right)  +\mu N_{2}\left(  s\right)  =0$, then
\[
\left\langle L_{1},N_{1}\right\rangle =1=\left\langle -\mu L_{2}%
,N_{1}\right\rangle =0
\]
which is a contradiction. Therefore the curve $\overline{\alpha}$ is a regular
curve. Then there exist a regular map $\varphi:s\rightarrow\overline{s}$
defined by%
\[
\overline{s}=\varphi\left(  s\right)  =%
{\displaystyle\int\limits_{0}^{s}}
\left\langle \frac{d^{3}\overline{\alpha}\left(  s\right)  }{ds^{3}}%
,\frac{d^{3}\overline{\alpha}\left(  s\right)  }{ds^{3}}\right\rangle
^{\frac{1}{6}}ds
\]

\noindent where $\overline{s}$ denotes the pseudo-arc parameter of
$\overline{\alpha}$, and we obtain%
\begin{equation}
\frac{d\overline{s}}{ds}=\frac{d\varphi\left(  s\right)  }{ds}=1. \label{a6}%
\end{equation}

\noindent So we can rewrite (\ref{a5}) as follows%
\begin{equation}
\overline{\alpha}\left(  \overline{s}\right)  =\alpha\left(  s\right)  +\mu
W_{3}\left(  s\right)  . \label{a7}%
\end{equation}
Differentiating (\ref{a7}) with respect to $s$, we get
\begin{equation}
\overline{L_{1}}\left(  \overline{s}\right)  =L_{1}\left(  s\right)  +\mu
N_{2}\left(  s\right)  . \label{a8}%
\end{equation}
If we differentiate (\ref{a8}), we get
\begin{equation}
\overline{L_{2}}\left(  \overline{s}\right)  =L_{2}\left(  s\right)  +\mu
N_{1}\left(  s\right)  . \label{a9}%
\end{equation}
Differentiating (\ref{a9}) once more, we arrive at%
\begin{equation}
\overline{W_{3}}\left(  \overline{s}\right)  =W_{3}\left(  s\right)  .
\label{a10}%
\end{equation}
This completes the proof.

\begin{example}
Let $\alpha$ be a null curve in $R_{2}^{5}$ defined by%
\[
\alpha\left(  s\right)  =\left(  \frac{s-s^{5}}{4\sqrt{15}},\frac{s^{2}+s^{4}%
}{4\sqrt{6}},\frac{s^{3}}{6},\frac{s^{2}-s^{4}}{4\sqrt{6}},\frac{s+s^{5}%
}{4\sqrt{15}}\right)  .
\]

\end{example}

\noindent Then we get the Cartan frame and the Cartan curvatures as follows:%
\begin{align*}
L_{1}  &  =\left(  \frac{1-5s^{4}}{4\sqrt{15}},\frac{2s+4s^{3}}{4\sqrt{6}%
},\frac{s^{2}}{2},\frac{2s-4s^{3}}{4\sqrt{6}},\frac{1+5s^{4}}{4\sqrt{15}%
}\right) \\
L_{2}  &  =\left(  \frac{-5s^{3}}{\sqrt{15}},\frac{2+12s^{2}}{4\sqrt{6}%
},s,\frac{2-12s^{2}}{4\sqrt{6}},\frac{5s^{3}}{\sqrt{15}}\right) \\
W_{3}  &  =\left(  -\sqrt{15}s^{2},\sqrt{6}s,1,-\sqrt{6}s,\sqrt{15}%
s^{2}\right) \\
N_{2}  &  =\left(  -2\sqrt{15}s,\sqrt{6},0,-\sqrt{6},2\sqrt{15}s\right) \\
N_{1}  &  =\left(  -2\sqrt{15},0,0,0,2\sqrt{15}\right) \\
k_{1}  &  =k_{2}=0.
\end{align*}

\noindent To get \ the Bertrand mate of $\alpha$ we can choose any nonzero
real number as $\mu$ so that its Bertrand mate is given by%
\begin{align*}
\overline{\alpha}\left(  \overline{s}\right)   &  =(\frac{\overline{s}%
-60\mu\left(  \overline{s}\right)  ^{2}-\left(  \overline{s}\right)  ^{5}%
}{4\sqrt{15}},\frac{24\mu\overline{s}+\left(  \overline{s}\right)
^{2}+\left(  \overline{s}\right)  ^{4}}{4\sqrt{6}},\frac{\left(  \overline
{s}\right)  ^{3}+6\mu}{6},\\
&  \frac{-24\mu\overline{s}+\left(  \overline{s}\right)  ^{2}-\left(
\overline{s}\right)  ^{4}}{4\sqrt{6}},\frac{\overline{s}+60\mu\left(
\overline{s}\right)  ^{2}+\left(  \overline{s}\right)  ^{5}}{4\sqrt{15}}),
\end{align*}

\noindent where $\overline{s}$ is the pseudo-arc parameter of $\overline
{\alpha}$, and a regular map $\varphi:s\rightarrow\overline{s}$ is given by%
\[
\overline{s}=\varphi\left(  s\right)  =s.
\]

\section{Pseudo-spherical Null Curves in $R_{2}^{n}$}

In this section, we characterize pseudo-spherical null curves in
$R_{2}^{n}$ in terms of the curvature functions (see \cite{sakaki}
for the pseudo-spherical curves in $R_{1}^{n}$). Here we only deal
with a family-type of null Cartan curves with the degeneration
degree two and nullity degree sequence is $\left\{
0,1,2,2,1,0,...,0\right\}  .$ The pseudo-sphere of
radius $r$ and center $\xi_{0}$ is given by%
\[
S_{2}^{n}\left(  r\right)  =\left\{  x\in R_{2}^{n}\mid\left\langle x-\xi
_{0},x-\xi_{0}\right\rangle =r^{2}\right\}
\]

For a Cartan curve $\alpha$ in $R_{2}^{n}$ parametrized by the pseudo-arc with
Cartan curvatures $\left\{  k_{1},k_{2},k_{3},...,k_{n-3}\right\}  $ and
$k_{n-3}\neq0,$ let us define a sequence of functions $\left\{  a_{1}%
,a_{2},...,a_{n-4}\right\}  $ inductively by%
\[
a_{1}=0,\text{ \ }a_{2}=\frac{1}{k_{3}},\text{ \ }a_{3}=-\frac{k_{3}^{\prime}%
}{\left(  k_{3}\right)  ^{2}k_{4}},\text{ \ }a_{i-1}=\frac{1}{k_{i}}\left(
a_{i-2}^{\prime}+a_{i-3}k_{i-1}\right)  ,\text{ \ }5\leqslant i\leqslant n-3.
\]

The Cartan equations for the curve mentioned above are as follows:%
\begin{align*}
\alpha^{\prime}  &  =L_{1}\\
L_{1}^{\prime}  &  =L_{2}\\
L_{2}^{\prime}  &  =W_{3}\\
W_{3}^{\prime}  &  =-k_{1}L_{2}+N_{2}\\
N_{2}^{\prime}  &  =k_{2}L_{1}+N_{1}-k_{1}W_{3}\\
N_{1}^{\prime}  &  =k_{2}L_{2}+k_{3}W_{4}\\
W_{4}^{\prime}  &  =-k_{3}L_{1}+k_{4}W_{5}\\
W_{i}^{\prime}  &  =-k_{i-1}W_{i-1}+k_{i}W_{i+1},\text{ \ \ \ }5\leqslant
i\leqslant n-3\\
W_{n-2}^{\prime}  &  =-k_{n-3}W_{n-3}.
\end{align*}

\noindent where $N_{1},N_{2}$ are null, $\left\langle L_{1},N_{1}\right\rangle
=1,$ $\left\langle L_{2},N_{2}\right\rangle =-1,$ $\left\{  L_{1},L_{2}%
,N_{1},N_{2}\right\}  $ and $\left\{  W_{3},W_{4},...,W_{n-3},W_{n-2}\right\}
$ are orthogonal, $\left\{  W_{3},W_{4},...,W_{n-3},W_{n-2}\right\}  $ is
orthonormal. We assume that $\left\{  \alpha^{\left(  i\right)  }\right\}
_{1\leqslant i\leqslant n}$ is positively oriented.

\begin{theorem}
Let $\alpha\left(  t\right)  $ be a null Cartan curve in $R_{2}^{n}$
parametrized with the pseudo-arc such that $k_{n-3}\neq0.$
\end{theorem}

\noindent\textbf{a)} If $\alpha\left(  t\right)  $ lies on a pseudo-sphere of
radius $r$, then
\[%
{\displaystyle\sum\limits_{i=2}^{n-4}}
a_{i}^{2}=r^{2}.
\]

\noindent\textbf{b)} If $a_{n-4}\neq0$ and $%
{\displaystyle\sum\limits_{i=2}^{n-4}}
a_{i}^{2}=r^{2}$ for some positive constant $r$, then $\alpha\left(  t\right)
$ lies on a pseudo-sphere of radius $r.$

\noindent\textbf{Proof.}

\noindent\textbf{a) }Suppose that $\alpha\left(  t\right)  $ lies on a
pseudo-sphere of radius $r$. Then there exist a fixed point $\xi_{0}$ $\in$
$R_{2}^{n}$ satisfying the following%
\begin{equation}
\left\langle \xi_{0}-\alpha\left(  t\right)  ,\text{ }\xi_{0}-\alpha\left(
t\right)  \right\rangle =r^{2}. \label{2}%
\end{equation}
Set%
\[
\xi_{0}-\alpha\left(  t\right)  =aL_{1}+bL_{2}+cN_{1}+dN_{2}+x_{1}W_{3}%
+x_{2}W_{4}+...+x_{n-4}W_{n-2}.
\]

\noindent Differentiating (\ref{2}), we have
\begin{equation}
\left\langle -L_{1},\xi_{0}-\alpha\left(  t\right)  \right\rangle =0,
\label{3}%
\end{equation}
and $c=0.$ Differentiating (\ref{3}), we have
\begin{equation}
\left\langle -L_{2},\xi_{0}-\alpha\left(  t\right)  \right\rangle =0 \label{4}%
\end{equation}
and $d=0$. Differentiating (\ref{4}), we have
\begin{equation}
\left\langle -W_{3},\xi_{0}-\alpha\left(  t\right)  \right\rangle =0,
\label{5}%
\end{equation}
and $a_{1}=x_{1}=0.$ Differentiating (\ref{5}), we have
\begin{equation}
\left\langle k_{1}L_{2}-N_{2},\xi_{0}-\alpha\left(  t\right)  \right\rangle =0
\label{6}%
\end{equation}
and $b=0.$ Differentiating (\ref{6}), we have
\begin{equation}
\left\langle k_{1}^{\prime}L_{2}+2k_{1}W_{3}-k_{2}L_{1}-N_{1},\xi_{0}%
-\alpha\left(  t\right)  \right\rangle =0 \label{7}%
\end{equation}
and $a=0.$ Differentiating (\ref{7}), we have%
\begin{equation}
\left\langle \left(  k_{1}^{\prime\prime}-2k_{1}^{2}-2k_{2}\right)
L_{2}+3k_{1}^{\prime}W_{3}+2k_{1}N_{2}-k_{2}^{\prime}L_{1}-k_{3}W_{4},\xi
_{0}-\alpha\left(  t\right)  \right\rangle =0 \label{8}%
\end{equation}
and $a_{2}=x_{2}=\frac{1}{k_{3}}.$ Differentiating (\ref{8}), we have
\begin{equation}
\left\langle ...-k_{3}^{\prime}W_{4}+k_{3}^{2}L_{1}-k_{3}k_{4}W_{5},\xi
_{0}-\alpha\left(  t\right)  \right\rangle =0 \label{9}%
\end{equation}
and $a_{3}=x_{3}=-\frac{k_{3}^{\prime}}{k_{3}^{2}k_{4}}.$ For $5\leqslant
i\leqslant n-3$, differentiating
\[
\left\langle \xi_{0}-\alpha\left(  t\right)  ,W_{i}\right\rangle =x_{i-2}%
\]

\noindent we have
\[
\left\langle -L_{1},W_{i}\right\rangle +\left\langle \xi_{0}-\alpha\left(
t\right)  ,-k_{i-1}W_{i-1}+k_{i}W_{i+1}\right\rangle =x_{i-2}^{\prime}%
\]
and%
\[
-x_{i-3}k_{i-1}+x_{i-1}k_{i}=x_{i-2}^{\prime}.
\]
So we get%
\[
x_{i-1}=\frac{1}{k_{i}}\left(  x_{i-2}^{\prime}+x_{i-3}k_{i-1}\right)  ,\text{
\ \ \ \ }5\leqslant i\leqslant n-3.
\]
We already have $x_{2}=a_{2},$ \ $x_{3}=a_{3}$. If we use the definition of
$\left\{  a_{i}\right\}  $, we can get $x_{i}=a_{i}$ for $1\leqslant
i\leqslant n-4.$ Therefore we have
\[
\xi_{0}-\alpha\left(  t\right)  =a_{2}W_{4}+a_{3}W_{5}+a_{4}W_{6}%
+...+a_{n-4}W_{n-2}%
\]
and by (\ref{2}), we have
\[%
{\displaystyle\sum\limits_{i=2}^{n-4}}
a_{i}^{2}=r^{2}.
\]
\textbf{b) }Suppose $a_{n-4}\neq0$ and
\begin{equation}%
{\displaystyle\sum\limits_{i=2}^{n-4}}
a_{i}^{2}=r^{2} \label{10}%
\end{equation}
for some positive constant $r$. Set
\[
\sigma\left(  t\right)  =\alpha\left(  t\right)  +a_{2}W_{4}+a_{3}%
W_{5}+...+a_{n-4}W_{n-2}%
\]
If we use the Frenet equations and the definition of $\left\{  a_{i}\right\}
$, we obtain%
\begin{align*}
\sigma^{\prime}\left(  t\right)   &  =\left(  1-a_{2}k_{3}\right)
L_{1}+\left(  a_{2}^{\prime}-a_{3}k_{4}\right)  W_{4}+\left(  a_{2}k_{4}%
+a_{3}^{\prime}-a_{4}k_{5}\right)  W_{5}\\
&  +\left(  a_{3}k_{5}+a_{4}^{\prime}-a_{5}k_{6}\right)  W_{6}+...+\left(
a_{n-6}k_{n-4}+a_{n-5}^{\prime}-a_{n-4}k_{n-3}\right)  W_{n-3}\\
&  +\left(  a_{n-5}k_{n-3}+a_{n-4}^{\prime}\right)  W_{n-2}\\
&  =\left(  a_{n-5}k_{n-3}+a_{n-4}^{\prime}\right)  W_{n-2}%
\end{align*}
Differentiating (\ref{10}), we have
\[
a_{2}a_{2}^{\prime}+a_{3}a_{3}^{\prime}+...+a_{n-5}a_{n-5}^{\prime}%
+a_{n-4}a_{n-4}^{\prime}=0.
\]
Using it together with the definition of $\left\{  a_{i}\right\}  ,$ we get%
\begin{align*}
a_{n-4}\left(  a_{n-4}^{\prime}+k_{n-3}a_{n-5}\right)   &  =k_{n-3}%
a_{n-4}a_{n-5}-a_{n-5}a_{n-5}^{\prime}\\
&  -a_{n-6}a_{n-6}^{\prime}-...-a_{4}a_{4}^{\prime}\\
&  -a_{3}a_{3}^{\prime}-a_{2}a_{2}^{\prime}\\
&  =a_{n-5}a_{n-6}k_{n-4}-a_{n-6}a_{n-6}^{\prime}-...\\
&  -a_{4}a_{4}^{\prime}-a_{3}a_{3}^{\prime}-a_{2}a_{2}^{\prime}\\
&  =a_{n-6}a_{n-7}k_{n-5}-a_{n-7}a_{n-7}^{\prime}-...\\
&  -a_{4}a_{4}^{\prime}-a_{3}a_{3}^{\prime}-a_{2}a_{2}^{\prime}\\
&  .\\
&  .\\
&  .\\
&  =a_{8}a_{7}k_{9}-a_{7}a_{7}^{\prime}-a_{6}a_{6}^{\prime}-a_{5}a_{5}%
^{\prime}-\\
&  -a_{4}a_{4}^{\prime}-a_{3}a_{3}^{\prime}-a_{2}a_{2}^{\prime}\\
&  .\\
&  .\\
&  .\\
&  =a_{4}a_{3}k_{5}-a_{3}a_{3}^{\prime}-a_{2}a_{2}^{\prime}\\
&  =a_{3}a_{2}k_{4}-a_{2}a_{2}^{\prime}\\
&  =a_{2}\left(  a_{3}k_{4}-a_{2}^{\prime}\right)  =0.
\end{align*}
So $\sigma^{\prime}\left(  t\right)  =0$ and $\sigma\left(  t\right)  =\xi
_{0}$ for some fixed point $\xi_{0}$ $\in$ $R_{2}^{n}.$ Thus we have
\[
\xi_{0}-\alpha\left(  t\right)  =\sum_{i=2}^{n-4}a_{i}W_{i+2},
\]
and by (\ref{10}), we have
\[
\left\langle \xi_{0}-\alpha\left(  t\right)  ,\xi_{0}-\alpha\left(  t\right)
\right\rangle =r^{2}.
\]
Hence $\alpha$ lies on a pseudo-sphere of radius $r.$

\section{Evolutes and Involutes in $R_{2}^{6}$}

Let $\alpha$ be a null Cartan curve in $R_{2}^{6}$ with degeneration
degree two and nullity degree sequence $\left\{
0,1,2,2,1,0,0\right\}  $. Then the
Frenet equations of $\alpha$ are as follows:%
\begin{align*}
\alpha^{\prime}  &  =L_{1},\text{ \ }L_{1}^{\prime}=L_{2},\text{ \ }%
L_{2}^{\prime}=W_{3},\text{ \ }W_{3}^{\prime}=-k_{1}L_{2}+N_{2},\text{ \ }\\
N_{2}^{\prime}  &  =k_{2}L_{1}+N_{1}-k_{1}W_{3},\text{ \ }N_{1}^{\prime}%
=k_{2}L_{2}+k_{3}W_{4},\\
W_{4}^{\prime}  &  =-k_{3}L_{1},
\end{align*}
where $L_{1},L_{2},N_{1},N_{2}$ are null, $\left\langle L_{1},N_{1}%
\right\rangle =1,$ $\left\langle L_{2},N_{2}\right\rangle =-1,$ $\left\{
L_{1},L_{2},N_{1},N_{2}\right\}  $ and $\left\{  W_{3},W_{4}\right\}  $ are
orthogonal, $\left\{  W_{3},W_{4}\right\}  $ is orthonormal. We assume that
$\left\{  \alpha^{\prime},\alpha^{\prime\prime},\alpha^{\left(  3\right)
},...,\alpha^{\left(  6\right)  }\right\}  $ is positively oriented, $\left\{
\alpha^{\prime},\alpha^{\prime\prime},\alpha^{\left(  3\right)  }%
,...,\alpha^{\left(  6\right)  }\right\}  $ and $\left\{  {}\right.
L_{1},L_{2},$

\noindent$W_{3},N_{2},N_{1},W_{4}\left.  {}\right\}  $ have the same
orientation. So we have $k_{3}>0$.

In the case $k_{3}\left(  t\right)  \neq0,$ let us define the evolute of
$\alpha$ by
\[
E\left(  t\right)  =\alpha\left(  t\right)  +\frac{1}{k_{3}\left(  t\right)
}W_{4}\left(  t\right)  ,
\]
which is the center of the osculating sphere at $\alpha\left(  t\right)  .$ On
the other hand, for a spacelike curve $c$ in $R_{2}^{6},$ we define the
involute of $c$ from a point $c\left(  t_{0}\right)  $ by%
\[
I\left(  t\right)  =c\left(  t\right)  -s\left(  t\right)  T\left(  t\right)
\text{,}%
\]
where $s\left(  t\right)  $ is the arc length of $c\left(  t\right)  $ from
$c\left(  t_{0}\right)  $ and $T\left(  t\right)  =\frac{c^{\prime}\left(
t\right)  }{\left\vert c^{\prime}\left(  t\right)  \right\vert }$ is the unit
tangent vector of c at $c\left(  t\right)  $ (see at \cite{sakaki} for the
evolute and involute of a curve in $R_{1}^{4}$)$.$

\begin{theorem}
Let $\alpha$ be a null Cartan curve in $R_{2}^{6}$ with the pseudo-arc
parameter $t$ such that $k_{3}\left(  t\right)  \neq0$ and $\left(  \frac
{1}{k_{3}\left(  t\right)  }\right)  ^{\prime}\neq0.$ Then the evolute $E$ of
$\alpha$ is a spacelike curve in $R_{2}^{6}$, and the involute $I_{E}$ of $E$
from some point coincides with $\alpha.$
\end{theorem}

\noindent\textbf{Proof.} By the Frenet equations for the Cartan curve $\alpha
$, the evolute of $\alpha$ satisfies
\begin{align*}
E^{\prime}\left(  t\right)   &  =L_{1}+\left(  \frac{1}{k_{3}}\right)
^{\prime}W_{4}+\frac{1}{k_{3}\left(  t\right)  }\left(  -k_{3}\left(
t\right)  L_{1}\right)  \\
&  =\left(  \frac{1}{k_{3}\left(  t\right)  }\right)  ^{\prime}W_{4}.
\end{align*}
So
\[
\left\langle E^{\prime},E^{\prime}\right\rangle =\left(  \left(  \frac
{1}{k_{3}\left(  t\right)  }\right)  ^{\prime}\right)  ^{2}>0,
\]
and $E$ is a spacelike curve. We only consider the case where $\left(
\frac{1}{k_{3}\left(  t\right)  }\right)  ^{\prime}>0,$ because the case where
$\left(  \frac{1}{k_{3}\left(  t\right)  }\right)  ^{\prime}<0$ is similar.
Then $E$ has unit tangent vector $T_{E}=W_{4},$ and the arc length
$s_{E}\left(  t\right)  $ of $E\left(  t\right)  $ from $E\left(
t_{0}\right)  $ is given by
\[
s_{E}\left(  t\right)  =%
{\displaystyle\int\limits_{t_{0}}^{t}}
\left\vert E^{\prime}\right\vert \text{ }dt=%
{\displaystyle\int\limits_{t_{0}}^{t}}
\left(  \frac{1}{k_{3}\left(  t\right)  }\right)  ^{\prime}dt=\frac{1}%
{k_{3}\left(  t\right)  }-\frac{1}{k_{3}\left(  t_{0}\right)  }.
\]
Therefore we have%
\[
\frac{1}{k_{3}\left(  t\right)  }=s_{E}\left(  t\right)  +\frac{1}%
{k_{3}\left(  t_{0}\right)  },
\]
which is the arc length of $E\left(  t\right)  $ from another point $E\left(
t_{1}\right)  .$

The involute $I_{E}$ of $E$ from $E\left(  t_{1}\right)  $ satisfies
\begin{align*}
I_{E}\left(  t\right)   &  =E\left(  t\right)  -\left(  s_{E}\left(  t\right)
+\frac{1}{k_{3}\left(  t_{0}\right)  }\right)  T_{E}\\
&  =\alpha\left(  t\right)  +\frac{1}{k_{3}\left(  t\right)  }W_{4}\left(
t\right)  -\frac{1}{k_{3}\left(  t\right)  }W_{4}\left(  t\right)
=\alpha\left(  t\right)  .
\end{align*}
Thus, we get the conclusion.

\begin{theorem}
Let $c$ be a spacelike curve in $R_{2}^{6}$ with pseudo-arc parameter $s$ such
that $c^{\prime\prime}\left(  s\right)  $ is null, $\left\langle c^{4}%
,c^{4}\right\rangle \neq0,$ and $\left\{  c^{\prime\prime},c^{\left(
3\right)  },...,c^{\left(  6\right)  }\right\}  $ is linearly independent.
Then for $s>0$ the involute $I$ of $c$ is a Cartan curve in $R_{2}^{6}$, and
the evolute $E_{I}$ of $I$ coincides with $c.$
\end{theorem}

\noindent\textbf{Proof}. As $T^{\prime}=c^{\prime\prime}$ is null, we may view
$T\left(  s\right)  $ as a null curve in $R_{2}^{6},$ and we have
$\left\langle T^{\left(  3\right)  },T^{\left(  3\right)  }\right\rangle
\geqslant0$ (look at \cite{intex-2})$.$ So the assumption $\left\langle
T^{\left(  3\right)  },T^{\left(  3\right)  }\right\rangle =\left\langle
C^{\left(  4\right)  },C^{\left(  4\right)  }\right\rangle \neq0$ implies that
$\left\langle T^{\left(  3\right)  },T^{\left(  3\right)  }\right\rangle >0.$
The involute $I\left(  s\right)  =$ $c\left(  s\right)  -sT\left(  s\right)  $
of the spacelike curve $c$ satisfies
\begin{align*}
I^{\prime}\left(  s\right)   &  =-sT^{\prime}\left(  s\right)  ,\text{
\ \ \ }I^{\prime\prime}\left(  s\right)  =-T^{\prime}\left(  s\right)
-sT^{\prime\prime}\left(  s\right)  ,\\
I^{\left(  3\right)  }\left(  s\right)   &  =-2T^{\prime\prime}\left(
s\right)  -sT^{\left(  3\right)  }\left(  s\right)  ,\text{ \ \ \ }I^{\left(
4\right)  }\left(  s\right)  =-3T^{\left(  3\right)  }\left(  s\right)
-sT^{\left(  4\right)  }\left(  s\right)  ,\\
I^{\left(  5\right)  }\left(  s\right)   &  =-4T^{\left(  4\right)  }\left(
s\right)  -sT^{\left(  5\right)  }\left(  s\right)  .
\end{align*}
For $s>0,$ $I$ is a null curve and
\[
\left\langle I^{\left(  3\right)  },I^{\left(  3\right)  }\right\rangle
=s^{2}\left\langle T^{\left(  3\right)  },T^{\left(  3\right)  }\right\rangle
>0.
\]
Let us denote $\left\langle T^{\left(  3\right)  },T^{\left(  3\right)
}\right\rangle ^{\frac{1}{2}}$ by $\eta.$ The pseudo-arc length $\nu\left(
s\right)  $ of $I$ is given by
\[
\nu\left(  s\right)  =%
{\displaystyle\int\limits_{s_{0}}^{s}}
\left\langle I^{\left(  3\right)  },I^{\left(  3\right)  }\right\rangle
^{\frac{1}{6}}ds=%
{\displaystyle\int\limits_{s_{0}}^{s}}
s^{\frac{1}{3}}\eta^{\frac{1}{3}}ds,
\]
and
\[
\frac{d\nu}{ds}=s^{\frac{1}{3}}\eta^{\frac{1}{3}}.
\]
Since $\left\{  T^{\prime},T^{\prime\prime},T^{\left(  3\right)  },T^{\left(
4\right)  },T^{\left(  5\right)  }\right\}  =\left\{  C^{\prime\prime
},C^{\left(  3\right)  },C^{\left(  4\right)  },C^{\left(  5\right)
},C^{\left(  6\right)  }\right\}  $ is linearly independent $\left\{
I^{\prime},I^{\prime\prime},I^{\left(  3\right)  },I^{\left(  4\right)
},I^{\left(  5\right)  }\right\}  $ is also linearly independent, and the null
curve $I$ is a Cartan curve with pseudo-arc length $\nu\left(  s\right)  $.
Let $\left\{  L_{1},L_{2},W_{3},N_{2},N_{1},W_{4}\right\}  $ be the Frenet
frame for the Cartan curve $I$ with Cartan curvatures $\left\{  k_{1}%
,k_{2},k_{3}\right\}  .$ Then we have
\begin{align*}
L_{1} &  =\frac{dI}{d\nu}=\frac{dI}{ds}\frac{ds}{d\nu}=-s^{\frac{2}{3}}%
\eta^{-\frac{1}{3}}T^{\prime},\\
L_{2} &  =\frac{dL_{1}}{d\nu}=\frac{dL_{1}}{ds}\frac{ds}{d\nu}=\frac{1}%
{3}\left(  s^{\frac{1}{3}}\eta^{-\frac{5}{3}}\eta^{\prime}-2s^{-\frac{2}{3}%
}\eta^{-\frac{2}{3}}\right)  T^{\prime}-s^{\frac{1}{3}}\eta^{-\frac{2}{3}%
}T^{\prime\prime},\\
W_{3} &  =\frac{dL_{2}}{d\nu}=\frac{1}{3}\left(  \frac{5}{3}s^{-1}\eta
^{-2}\eta^{\prime}-\frac{5}{3}\eta^{-3}\left(  \eta^{\prime}\right)  ^{2}%
+\eta^{-2}\eta^{\prime\prime}+\frac{4}{3}s^{-2}\eta^{-1}\right)  T^{\prime}\\
&  +\left(  \eta^{-2}\eta^{\prime}-s^{-1}\eta^{-1}\right)  T^{\prime\prime
}-\eta^{-1}T^{\left(  3\right)  },\\
N_{2} &  =\frac{dW_{3}}{d\nu}+k_{1}L_{2},\\
N_{1} &  =\frac{dN_{2}}{d\nu}-k_{2}L_{1}+k_{1}W_{3}.
\end{align*}
Note that
\[
\text{\emph{span}}\left\{  L_{1},L_{2},W_{3},N_{2},N_{1}\right\}
=\text{\emph{span}}\left\{  T^{\prime},T^{\prime\prime},T^{\left(  3\right)
},T^{\left(  4\right)  },T^{\left(  5\right)  }\right\}
\]
and
\[
\left\langle T,T^{\prime}\right\rangle =\left\langle T,T^{\prime\prime
}\right\rangle =\left\langle T,T^{\left(  3\right)  }\right\rangle
=\left\langle T,T^{\left(  4\right)  }\right\rangle =\left\langle T,T^{\left(
5\right)  }\right\rangle =0
\]
We obtain $W_{4}=\mp T$. We can take $W_{4}=T,$ because the case $W_{4}=-T$ is
similar. Then
\begin{align*}
k_{3} &  =-\left\langle \frac{dW_{4}}{d\nu},N_{1}\right\rangle =-\left\langle
\frac{dT}{d\nu},N_{1}\right\rangle \\
&  =-\left\langle s^{-\frac{1}{3}}\eta^{-\frac{1}{3}}T^{\prime},N_{1}%
\right\rangle \\
&  =\left\langle s^{-1}L_{1},N_{1}\right\rangle =\frac{1}{s}.
\end{align*}
Thus the evolute $E_{I}$ of $I$ satisfies
\[
E_{I}\left(  s\right)  =I\left(  s\right)  +\frac{1}{k_{3}\left(  s\right)
}W_{4}=c\left(  s\right)  -sT+sT=c\left(  s\right)  ,
\]
which completes the proof.

\bigskip

\bigskip

\bigskip

\bigskip

\bigskip

\bigskip

\noindent{Mehmet G\"{o}\c{c}men}

\noindent{Department of Mathematics, Faculty of Arts and Sciences,
\.{I}n\"{o}n\"{u} University }

\noindent{44280 Malatya, Turkey }

\noindent{Email: mgocmen1903@gmail.com}

\noindent{Sad\i k Kele\c{s} }

\noindent{Department of Mathematics, Faculty of Arts and Sciences,
\.{I}n\"{o}n\"{u} University }

\noindent{44280 Malatya, Turkey }

\noindent{Email: skeles@inonu.edu.tr }

\end{document}